\documentclass[12pt]{amsart}

\newcommand{\norm}[1]{\left\Vert#1\right\Vert}
\newcommand{\abs}[1]{\left\vert#1\right\vert}

\newcommand{\set}[1]{\left\{#1\right\}}

\newcommand{\Part}{\mathcal{P}}
\newcommand{\state}[1]{\varphi[#1]}
\newcommand{\Prod}{\mathrm{Pr}}
\newcommand{\St}{\mathrm{St}}
\newcommand{\NC}{\mathit{NC}}

\newcommand{\mf}[1]{\mathbb{#1}}
\newcommand{\mb}[1]{\mathbf{#1}}
\newcommand{\mc}[1]{\mathcal{#1}}
\DeclareMathOperator{\Mob}{M\ddot{o}b}

\title[Free stochastic measures II]{Free stochastic measures via noncrossing partitions II}

\author[M. Anshelevich]{Michael Anshelevich}
\address{Department of Mathematics, University at California, Berkeley, CA 94720}
\email{manshel@math.berkeley.edu}

\subjclass{Primary 46L54; Secondary 60G, 60H}

\newtheorem*{Main}{Main Theorem}
\newtheorem{Thm}{Theorem}
\newtheorem{Cor}[Thm]{Corollary}
\newtheorem{Lemma}[Thm]{Lemma}

\theoremstyle{definition}
\newtheorem{Defn}{Definition}
\newtheorem{Notation}[Defn]{Notation}

\theoremstyle{remark}
\newtheorem{Remark}{Remark}
\newtheorem{Ex}[Remark]{Example}
  
\date{August 7, 2001}

\begin{document}

\begin{abstract} 
We show that for stochastic processes with freely independent increments, the partition-dependent stochastic measures can be expressed purely in terms of the higher stochastic measures and the higher diagonal measures of the original process. 
\end{abstract}

\maketitle

\section{Introduction}

Starting with an operator-valued stochastic process with freely independent increments $X(t)$, in \cite{Ans00} we defined two families $\set{\Prod_\pi}$ and $\set{\St_\pi}$ indexed by set partitions. These objects give a precise meaning to the following heuristic expressions. For a partition $\pi = (B_1, B_2, \ldots, B_n) \in \Part(k)$, temporarily denote by $c(i)$ the number of the class $B_{c(i)}$ to which $i$ belongs. Then, heuristically,
\begin{equation*}
\Prod_\pi(t) = \int_{[0, t)^n} dX(s_{c(1)}) dX(s_{c(2)}) \cdots dX(s_{c(k)})
\end{equation*}
and
\begin{equation*}
\St_\pi(t) = \int_{\substack{[0, t)^n \\ \text{all $s_i$'s distinct}}} dX(s_{c(1)}) dX(s_{c(2)}) \cdots dX(s_{c(k)}).
\end{equation*}
In particular, denote by $\psi_k$ and $\Delta_k$ the higher stochastic measures and the higher diagonal measures, defined, respectively, by 
\begin{equation*}
\psi_k(t) = \int_{\substack{[0, t)^k \\ \text{all $s_i$'s distinct}}} dX(s_1) dX(s_2) \cdots dX(s_k)
\end{equation*}
and 
\begin{equation*}
\Delta_k(t) = \int_{[0, t)} (dX(s))^k.
\end{equation*}
Rigorous definitions of all these objects in terms of Riemann sums are given below. These definitions were motivated by \cite{RW97}, where corresponding objects were defined for the usual L\'{e}vy processes. There is a number of differences between the classical and the free case. First, the free increments property implies that $\St_\pi = 0$ unless $\pi$ is a noncrossing partition. Second, the point of the analysis of \cite{RW97} was that while we are really interested in the stochastic measures $\St_\pi$, notably $\psi_k$, these are rather hard to define or to handle. However, by the use of M\"{o}bius inversion these can be expressed through the $\Prod_\pi$. It is easy to see that if the increments of the process $X$ commute, in the defining expression for $\Prod_\pi$ all the terms corresponding to the same class can be collected together, and the result is just a product measure over the classes of the partition, $\Prod_\pi = \Delta_{B_1} \Delta_{B_2} \cdots \Delta_{B_n}$. So in this way stochastic measures $\St$ can be defined using ordinary product measures. This fact is a consequence of the commutativity of the increments of the process; in the free probability case the operators do not commute, and unless the classes of the partition $\pi$ are just intervals we cannot expect $\Prod_\pi$ to be a product measure; indeed a counterexample was given in \cite{Ans00}.

In this paper we show that while we cannot expect nice factorization properties in the general case of noncommuting variables, the free independence of the increments does imply a product-like property. Namely, by an argument similar to the above one, if the increments of the process commute, then $\St_\pi = \psi_k(\Delta_{B_1}, \Delta_{B_2}, \ldots, \Delta_{B_n})$. This property certainly does not hold either if the increments do not commute, but if the increments are freely independent it can be modified as follows. In a noncrossing partition, one distinguishes classes which are \emph{inner}, or covered by some other classes, and \emph{outer}. For example, in the noncrossing partition $((1, 6, 7) (2, 5) (3) (4) (8) (9, 10))$, the classes $\set{(2,5), (3), (4)}$ are inner while the classes $\set{(1, 6, 7), (8), (9, 10)}$ are outer. For a partition with only outer classes, which therefore have to be intervals, the product decomposition of $\Prod_\pi$ and the above decomposition of $\St_\pi$ hold even in the noncommutative case. We show in the main theorem of this paper that the inner classes, while making a complicated contribution to $\Prod_\pi$, make only scalar contributions to $\St_\pi$, and those contributions commute with everything.

We also use the opportunity to make extensions of the definitions of \cite{Ans00} in various directions. The objects $\St_\pi$ and $\Prod_\pi$ were defined by limits of Riemann-like sums with uniform subdivisions. In this paper we extend that definition to arbitrary subdivisions. We also make some preliminary steps towards defining multi-dimensional free stochastic measures.

\noindent\textbf{Acknowledgments.} This paper stems from the work started during the special semester ``Free probability and operator spaces'' at the Institut Henri Poincar\'{e} -- Centre Emile Borel, with the support from the Fannie and John Hertz Foundation; it was completed with the support from an NSF postdoctoral fellowship. The author is grateful to the referee for numerous comments and suggestions.

\section{Preliminaries}

This paper is a sequel to \cite{Ans00}; see that paper for all the definitions that are not explicitly provided here.

\subsection{Notation}

Denote by $[n]$ the set $\set{1, 2, \ldots, n}$.

For two vectors $\mb{X} = (X_1, X_2, \ldots, X_k)$ and $\mb{Z} = (Z_1, Z_2, \ldots, Z_{k-1})$ denote 
\begin{equation*}
X \circ Z = X_1 Z_1 X_2 Z_2 \cdots X_{k-1} Z_{k-1} X_k.
\end{equation*}

For a collection of vectors $\set{\mb{X}_i}_{i=1}^n$, denote by $(\mb{X}_1, \mb{X}_2, \ldots, \mb{X}_n)$ their concatenation.

For a collection of objects $\set{y^{(i)}_j}$ and a multi-index $\bar{v} = (v_1, v_2, \ldots, v_n)$, we will throughout the paper use the notation $\mb{y}_{\bar{v}}$ to denote $\prod_{j=1}^n y_{v_j}^{(j)}$.

For a family of functions $\set{F_j}$, where $F_j$ is a function of $j$ arguments, $\bar{v}$ a vector with $k$ components, and $B \subset [k]$, denote $F(\bar{v}) = F_k(\bar{v})$ and 
\begin{equation*}
F(B; \bar{v}) = F_{\abs{B}}(v_{i(1)}, v_{i(2)}, \ldots, v_{i(\abs{B})}),
\end{equation*}
where $B = (i(1), i(2), \ldots, i(\abs{B}))$. In particular, using this notation $\mb{y}_{(B; \bar{v})} = \prod_{i \in B} y^{(i)}_{v_i}$.

\subsection{Partitions}

Denote by $\Part(k)$ and $\NC(k)$ the lattice of all set partitions of the set $[k]$ and its sub-lattice of noncrossing partitions. Let $\hat{0}$ and $\hat{1}$ be the smallest and the largest elements in the lattice ordering, namely $\hat{0} = ((1), (2), \ldots, (k))$ and $\hat{1} = (1, 2, \ldots, k)$. Denote by $\wedge$ the join operation in the lattices. For $\pi \in \NC(k)$, denote by $K(\pi)$ its Kreweras complement. For $\pi \in \Part(n)$, define $\pi^{op} \in \Part(n)$ to be $\pi$ taken in the opposite order, i.e.
\begin{equation*}
i  \stackrel{\pi^{op}}{\sim} j \; \Leftrightarrow \; (n-i+1) \stackrel{\pi}{\sim} (n-j+1).
\end{equation*}
For $\pi \in \Part(n), \sigma \in \Part(k)$, define $\pi + \sigma \in \Part(n+k)$ by
\begin{equation*}
i \stackrel{\pi + \sigma}{\sim} j \Leftrightarrow ((i, j \leq n, i \stackrel{\pi}{\sim} j)
\text{ or } (i, j > n, (i-n) \stackrel{\sigma}{\sim} (j-n))).
\end{equation*}

\subsection{Free cumulants}

All the operators involved will live in an ambient noncommutative probability space $(\mc{A}, \varphi)$, where $\mc{A}$ is a finite von Neumann algebra, and $\varphi$ is a faithful normal tracial state on it. Let $\mb{A} = (A_1, A_2, \ldots, A_k)$ be a $k$-tuple of self-adjoint operators in $(\mc{A}, \varphi)$. Denote their joint moments by 
\begin{equation*}
M(\mb{A}) = \state{A_1 A_2 \cdots A_k}.
\end{equation*}
For a noncrossing partition $\pi$, denote using the above notation
\begin{equation*}
M_\pi(\mb{A}) = \prod_{B \in \pi} M(B; \mb{A}).
\end{equation*}
Also define the combinatorial $R$-transform, or the collection of joint free cumulants $R(\mb{A})$: denoting
\begin{equation*}
R_\pi(\mb{A}) = \prod_{B \in \pi} R(B; \mb{A}),
\end{equation*}
the functional $R$ is determined inductively by
\begin{equation*}
M(\mb{A}) = \sum_{\sigma \in \NC(k)} R_\sigma(\mb{A}),
\end{equation*}
or more generally by
\begin{equation*}
M_\pi(\mb{A}) = \sum_{\substack{\sigma \in \NC(k) \\ \sigma \leq \pi}} R_\sigma(\mb{A}).
\end{equation*}
Any such relation can be inverted by using M\"{o}bius inversion, so we also get
\begin{equation*}
R_\pi(\mb{A}) = \sum_{\substack{\sigma \in \NC(k) \\ \sigma \leq \pi}} \Mob(\sigma, \pi) M_\sigma(\mb{A}),
\end{equation*}
where $\Mob$ is the relative M\"{o}bius function on the lattice of noncrossing partitions. In particular, since $\abs{M_\pi(\mb{A})} \leq \prod_{i=1}^k \norm{A_i}$ and $\Mob(\pi, \sigma), \abs{\NC(k)}$ are products of Catalan numbers and so are bounded in norm by $4^k$, we conclude that $\abs{R_\pi(\mb{A})} \leq 16^k \prod_{i=1}^k \norm{A_i}$.

Finally, the relation between the free cumulants and free independence is expressed in the ``mixed cumulants are zero'' condition: $R(\mb{A}) = 0$ whenever some $A_i, A_j$ are freely independent.

\begin{Defn}
Let $\mb{X} = (X^{(1)}, X^{(2)}, \ldots, X^{(k)})$ be a $k$-tuple of free stochastic measures with distributions $\mu^{(1)}, \mu^{(2)}, \ldots, \mu^{(k)}$. Here, $\mu^{(i)}$ is a freely infinitely divisible distribution with compact support, and $X^{(i)}$ is an operator-valued measure on $\mf{R}$ that is self-adjoint, additive, stationary, and has freely independent increments. We say that the $k$-tuple is \emph{consistent} if the following extra conditions are satisfied.
\begin{enumerate}
\item
Free increments: For a family of disjoint intervals $\set{I_i}_{i=1}^n$ and a multi-index $\bar{u}$ on length $n$, the family $\set{X^{(u_i)}(I_i)}_{i=1}^n$ is a freely independent family.
\item 
Stationarity: For an interval $I$ and multi-index $\bar{u}$,
\[
\state{X^{(u_1)}(I) X^{(u_2)}(I) \ldots X^{(u_n)}(I)}
\]
depends only on $\bar{u}$ and $\abs{I}$.
\item
Continuity: For a fixed $\bar{u}$, the function
\[
\abs{I} \mapsto \state{X^{(u_1)}(I) X^{(u_2)}(I) \ldots X^{(u_n)}(I)}
\]
is continuous.
\end{enumerate}
\end{Defn}

Note that these conditions together imply that for an arbitrary collection of intervals $\set{I_i}_{i=1}^n$, $\state{X^{(u_1)}(I_1) X^{(u_2)}(I_2) \ldots X^{(u_n)}(I_n)}$ depends only on $\bar{u}$ and the sizes of all elements of $\set{\bigcap_{i \in G} I_i: G \subset [n]}$. Note also that by definition and M\"{o}bius inversion, the stationarity and continuity properties apply not just to $M$, but to $M_\pi, R, R_\pi$ as well.

\begin{Remark}
It is easy to see that a $k$-tuple $\mb{X}$ is consistent if $X^{(i)} = X$ for all $i$, or if the family $\set{X^{(i)}}_{i=1}^k$ is a freely independent family. More examples are given in Lemma \ref{Lemma:Diagonal} and in Remark \ref{Remark:hda}.
\end{Remark}

Fix a consistent $k$-tuple of free stochastic measures $\mb{X}$. For $t>0$, denote $X^{(i)}(t) = X^{(i)}([0,t))$. Throughout most of the paper, we will consider $t=1$, in which case we will omit $t$ from the notation; in particular we will denote $X^{(i)}([0,1))$ simply by $X^{(i)}$. Let $\mc{S} = (I_1, I_2, \ldots, I_N)$ be a subdivision of the interval $[0, t)$ into $N = \abs{\mc{S}}$ disjoint half-open intervals, listed in increasing order. Denote $\delta(\mc{S}) = \max_{1 \leq i \leq N} \abs{I_i}$. Let $X^{(i)}_j (\mc{S}) = X^{(i)}(I_j)$. In the future we will frequently omit the dependence on $\mc{S}$ and $N$ in the notation.

\begin{Notation}
For any set $G$ and a partition $\pi \in \Part(k)$, denote
\[
G^k_\pi = \set{\bar{v} \in G^k: i \stackrel{\pi}{\sim} j \Leftrightarrow v_i = v_j}
\]
and
\[
G^k_{\geq \pi} = \set{\bar{v} \in G^k: i \stackrel{\pi}{\sim} j \Rightarrow v_i = v_j}.
\] 
Denote 
\begin{equation*}
\St_{\pi}(\mb{X}, \mc{S}) = \sum_{\bar{v} \in [N]^k_{\pi}} \mb{X}_{\bar{v}} (\mc{S})
\end{equation*}
and
\begin{equation*}
\Prod_{\pi}(\mb{X}, \mc{S}) = \sum_{\bar{v} \in [N]^k_{\geq \pi}} \mb{X}_{\bar{v}} (\mc{S}).
\end{equation*}
\end{Notation}

\begin{Defn}
\label{Defn:fsm}
Define the \emph{free stochastic} and \emph{product measures} depending on a partition to be the limits along the net of subdivisions of the interval $[0, t)$
\begin{align*}
\St_{\pi}(\mb{X}, t) &= \lim_{\delta(\mc{S}) \rightarrow 0} \St_{\pi}(\mb{X}, \mc{S}), \\
\Prod_{\pi}(\mb{X}, t) &= \lim_{\delta(\mc{S}) \rightarrow 0} \Prod_{\pi}(\mb{X}, \mc{S}).
\end{align*}
In particular, let the higher diagonal measure be
\[
\Delta(\mb{X}, t) = \St_{\hat{1}}(\mb{X}, t) = \Prod_{\hat{1}}(\mb{X}, t),
\]
and the $k$-dimensional free stochastic measure be $\psi(\mb{X}, t) = \St_{\hat{0}}(\mb{X}, t)$. If $X^{(i)} = X$ for all $i$, we denote $\Delta(\mb{X}, t)$ by $\Delta_k(t)$ and $\psi(\mb{X}, t)$ by $\psi_k(t)$. 
\end{Defn}

Here the limits are taken in the operator norm; the proof of their existence is part of the arguments in the next section.

\begin{Lemma}
\label{lem:r}
For an arbitrary family of intervals $\set{J_i \subset [0, 1)}_{i=1}^n$, a multi-index $\bar{u}$, and $\pi \in \NC(k)$, 
\begin{equation*}
R_\pi(X^{(u_1)}(J_1), X^{(u_2)}(J_2), \ldots, X^{(u_n)}(J_n))
= \left( \prod_{B \in \pi} \abs{\bigcap_{i \in B} J_i} \right) R_\pi(X^{(u_1)}, X^{(u_2)}, \ldots, X^{(u_n)}).
\end{equation*}
In particular, for a subdivision $\mc{S} = (I_1, I_2, \ldots, I_N)$ of $[0, 1)$,
\begin{equation*}
R_\pi(X^{(u_1)}_j, X^{(u_2)}_j, \ldots, X^{(u_n)}_j) = \abs{I_j}^{\abs{\pi}} R_\pi(X^{(u_1)}, X^{(u_2)}, \ldots, X^{(u_n)}).
\end{equation*}
\end{Lemma}

\begin{proof}
The second statement follows from the first one with all intervals $J_i = I_j$ with a fixed $j$. For the first statement, it suffices to show that
\begin{equation*}
R(X^{(u_1)}(J_1), X^{(u_2)}(J_2), \ldots, X^{(u_n)}(J_n))
= \abs{\bigcap_{i = 1}^n J_i} R(X^{(u_1)}, X^{(u_2)}, \ldots, X^{(u_n)}).
\end{equation*}
Moreover, since each $X^{(j)}$ is an additive process with freely independent increments, the expression $R_\pi(A_1, A_2, \ldots, A_n)$ is multi-linear in its arguments, and all mixed cumulants are equal to $0$, it suffices to show that
\begin{equation}
\label{Rsmall}
R(X^{(u_1)}(I), X^{(u_2)}(I), \ldots, X^{(u_n)}(I)) 
= \abs{I} R(X^{(u_1)}, X^{(u_2)}, \ldots, X^{(u_n)})
\end{equation}
with $I = \bigcap_{i = 1}^n J_i$. First suppose that $I = I_j$, one of the intervals in a uniform subdivision, with $I_i = [\frac{i-1}{N}, \frac{i}{N})$. Then using the same properties as above,
\begin{align*}
R(X^{(u_1)}, X^{(u_2)}, \ldots, X^{(u_n)}) 
&= \sum_{\bar{v} \in [N]^n} R(X^{(u_1)}_{v_1}, X^{(u_2)}_{v_2}, \ldots, X^{(u_n)}_{v_n}) \\
&= \sum_{i=1}^N R(X^{(u_1)}_i, X^{(u_2)}_i, \ldots, X^{(u_n)}_i) \\
&= N R(X^{(u_1)}_j, X^{(u_2)}_j, \ldots, X^{(u_n)}_j),
\end{align*}
where in the last equality we have used that fact that by stationarity of $\mb{X}$, $R(X^{(u_1)}_i, \ldots, X^{(u_n)}_i)$ does not depend on $i$. Therefore
\begin{align*}
R(X^{(u_1)}_j, X^{(u_2)}_j, \ldots, X^{(u_n)}_j) 
&= N^{-1} R(X^{(u_1)}, X^{(u_2)}, \ldots, X^{(u_n)}) \\
&= \abs{I_j} R(X^{(u_1)}, X^{(u_2)}, \ldots, X^{(u_n)}).
\end{align*}
By stationarity, it follows that equation \eqref{Rsmall} holds for $\abs{I} = 1/N$ and consequently for any rational $\abs{I}$. The result follows for general $I$ by the continuity assumption on $\mb{X}$.
\end{proof}

This was the main fact used in the proofs of \cite{Ans00}, and so the results from that paper carry over to the consistent $k$-tuples of free stochastic measures. Note that all of these results were proven for uniform subdivisions. However, their proofs carry over to the more general definitions of this paper without difficulty; see Lemma \ref{lem:proj} for an example of a computation. We list some of those results, with the numbering from \cite{Ans00}.
\begin{enumerate}
\item
$\St_\pi(\mb{X}) =0$ unless $\pi$ is noncrossing (Theorem 1).
\item
\label{num1}
$\state{\St_\pi(\mb{X})} = R_\pi(\mb{X})$ (Corollary 2). 
\item
\label{num2}
If $\pi$ contains an inner singleton and $\mb{X}$ is centered, then $\St_\pi(\mb{X}) = 0$ (Proposition 1). 
\item
\label{num3}
If $Z$ is freely independent from the free stochastic measure $X$, then
\[
\lim_{\delta(\mc{S}) \rightarrow 0} \sum_{i=1}^N X_i Z X_i = \state{Z} \Delta_2(X)
\]
(Corollary 13). 
\item
The limit defining $\St_\pi$ exists in the norm topology if the corresponding limit exists for the free Poisson process. In particular, for any consistent $k$-tuple $\mb{X}$ of free stochastic measures, $\Delta(\mb{X})$ is well-defined. In fact, the argument in \cite{Ans00} needs to be modified ($\Prod$ should be used in place of $\St$); such a modification is contained in Lemma \ref{lem:limit} of this paper. 
\end{enumerate}

Results \ref{num1} and \ref{num2} are consequences of, and result \ref{num3} is parallel to, the following theorem.

\begin{Main}
\label{thm:main}
Let $\pi$ be a noncrossing partition of $[k]$ with $o(\pi)$ outer classes $B_1, \ldots, B_{o(\pi)}$ and $i(\pi)$ inner classes $C_1, \ldots, C_{i(\pi)}$. Then
\begin{equation*}
\St_{\pi}(\mb{X}) = \prod_{i=1}^{i(\pi)} R(C_i; \mb{X}) \cdot \psi(\Delta(B_1; \mb{X}), \Delta(B_2; \mb{X}), \ldots, \Delta(B_{o(\pi)}; \mb{X})).
\end{equation*}
\end{Main}

\begin{Remark}
The distinction between the inner and the outer classes of a noncrossing partition was noted in \cite{BLS96}. It would be interesting to see what the relation is between the conditionally free cumulants of that paper (which are scalar-valued) and our $\St_\pi$; cf.\ also \cite{Mlo99}.
\end{Remark}

\begin{Ex}
Let $\pi$ be as in the theorem, and $\set{X(t)}$ be the free Poisson process. Then for $n \geq 1$, $r_n(t) = t$, $\Delta_n(t) = X(t)$. Therefore the theorem states that
\begin{equation*}
\St_\pi(t) = t^{i(\pi)} \psi_{o(\pi)}(t).
\end{equation*}
\end{Ex}

\begin{Ex}
Let $\pi$ be as in the theorem, and $\set{X(t)}$ be the free Brownian motion. Then the free cumulants of $\mu_t$ are $r_1(t) = 0$, $r_2(t) = t$, $r_n(t) = 0$ for $n>2$, and the diagonal measures of $X$ are $\Delta_1(t)  = X(t)$, $\Delta_2(t) = t$, $\Delta_n(t) = 0$ for $n > 2$. Therefore the theorem states that 
\begin{align*}
\St_\pi(X, t) 
&= \begin{cases}
0 \text{ if $\pi$ contains a class of more than $2$ elements}, \\
0 \text{ if $\pi$ contains an inner singleton}, \\
t^{i(\pi)} t^{\#\set{B_j: \abs{B_j} = 2}} \psi_{\#\set{B_j: \abs{B_j} = 1}}(t)   \text{ otherwise}
\end{cases} \\
&= \begin{cases}
t^{\#\set{V \in \pi: \abs{V} = 2}} \psi_{\#\set{V \in \pi: \abs{V} = 1}}(t) \text{ if } \forall i, j, \abs{C_i} = 2, \abs{B_j} = 1, 2, \\
0 \text{ otherwise}.
\end{cases}
\end{align*} 
\end{Ex}

\section{Proof of the theorem}

We start the analysis with a single free Poisson stochastic measure. It has the following remarkable representation: one can take $I \mapsto X(I)$ to be $s p(I) s$, where $s$ is a variable with a semicircular distribution freely independent from $p(I)$, and $I \mapsto p(I)$ is a projection-valued measure, so that disjoint intervals correspond to orthogonal projections, and $\state{p(I)} = \abs{I}$.

\begin{Lemma}
\label{lem:proj}
Given a subdivision $\mc{S} = (I_1, I_2, \ldots, I_N)$ of $[0,1)$, let $\set{p_i}_{i=1}^N$ be orthogonal projections adding up to $1$ with $\state{p_i} = \abs{I_i}$. Let $\set{Z_{i,j}}_{i \in [N], j \in [k]}$ be a family of operators (dependent on $\mc{S}$) such that for each $i$, the family $\set{Z_{i,j}}_{j \in [k]}$ is freely independent from $p_i$. In addition, assume that for each $i$, at least one of $Z_{i,j}$ is centered, and that for all $i, j, \mc{S}$, $\norm{Z_{i,j}} < c/16$. Then
\begin{equation*}
 \lim_{\delta(\mc{S}) \rightarrow 0} \sum_{i=1}^N p_i Z_{i,1} p_i Z_{i,2} \cdots Z_{i,k} p_i = 0,
\end{equation*}
where the limit is taken in the operator norm.
\end{Lemma}

\begin{proof}
Denote 
\begin{equation*}
A(\mc{S}) = \sum_{i=1}^N p_i Z_{i,1} p_i Z_{i,2} \cdots Z_{i,k} p_i = \sum_{i=1}^N (p_i, \ldots, p_i) \circ \mb{Z}_i,
\end{equation*}
where $\mb{Z}_i = (Z_{i,1}, \ldots, Z_{i,k})$. Then 
\begin{equation*}
(A(\mc{S}) A(\mc{S})^\ast)^n = \sum_{i=1}^N (p_i, \ldots, p_i) \circ (\mb{Z}_i, (\mb{Z}_i)^\ast, \ldots, \mb{Z}_i, (\mb{Z}_i)^\ast),
\end{equation*} 
where $(\mb{Z}_i)^\ast = (Z_{i,k}^\ast, \ldots, Z_{i,1}^\ast)$. Then as in \cite[Theorem 3]{Ans00},
\begin{equation*}
\state{(A(\mc{S}) A(\mc{S})^\ast)^n}
= \sum_{i=1}^N \sum_{\pi \in \NC(2nk)} R_{K(\pi)}(\mb{Z}_i, (\mb{Z}_i)^\ast, \ldots, \mb{Z}_i, (\mb{Z}_i)^\ast)) \cdot M_\pi(p_i, p_i, \ldots, p_i).
\end{equation*}
At least one $Z_{i,j}$ is centered, so $\abs{K(\pi)} \leq 2n k - n$ and $\abs{\pi} \geq (n+1)$ (since $\abs{\pi} + \abs{K(\pi)} = 2nk +1$). Thus 
\[
\state{(A(\mc{S}) A(\mc{S})^\ast)^n}
\leq \sum_{i=1}^N \sum_{\substack{\pi \in \NC(2nk) \\ \abs{\pi} \geq n+1}} c^{2nk} \state{p_i}^{\abs{\pi}}
\leq 4^{2n k} c^{2n k} \delta(\mc{S})^n.
\]
Therefore
\begin{equation*}
\norm{A(\mc{S})}_{2n} < \delta(\mc{S}) ^{1/2k} 4^k c^k,
\end{equation*}
and so $\norm{A(\mc{S})}  < \delta(\mc{S})^{1/2k}  (4c)^k$, which converges to $0$ as $\delta(\mc{S}) \rightarrow 0$.
\end{proof}

\begin{Lemma}
\label{Lem:PoissonCentered}
Let $I \mapsto X(I) = s p(I) s$ be a free Poisson stochastic measure (see \cite{Ans00}), so that $X_i = s p_i s$ with $p_i$ as in the previous lemma. Let $Z_{i,j}$ be as in the previous lemma. Then 
\begin{equation*}
\lim_{\delta(\mc{S}) \rightarrow 0} \sum_{i=1}^N X_i Z_{i,1} X_i Z_{i,2} \ldots Z_{i,k} X_i = 0.
\end{equation*}
\end{Lemma}

\begin{proof}
By the free independence assumption on $\set{Z_{i,j}}$, the joint distribution of the $(k+1)$-tuple $\set{X_i, Z_{i,1}, Z_{i,2}, \ldots, Z_{i,k}}$ is, for each $i$, entirely determined by the distribution of $X_i$ and the joint distribution of $\set{Z_{i,j}}_{j=1}^k$. These distributions, and hence the conclusion of the lemma, are not changed if we assume in addition that the family $\set{Z_{i,j}}_{j=1}^k$ is freely independent from $s$.

$(X_i, \ldots, X_i) \circ \mb{Z}_i = (s p_i s, \ldots, s p_i s) \circ \mb{Z}_i = s ((p_i, \ldots, p_i) \circ (s \mb{Z}_i s)) s$, where $s \mb{Z}_i s$ denotes the vector $\mb{Z}_i$ with each term multiplied by $s$ on both sides. Since $s$ and $\mb{Z}_i$ are freely independent, if $Z_{i,j}$ is centered then so is $sZ_{i,j}s$. Then the previous lemma implies the result.
\end{proof}

\begin{Lemma}
\label{lem:st1}
Let $\pi \in \NC(k)$ have only one outer class $B$ consisting of $n+1$ elements. That is,  
\begin{equation*}
\pi = \set{(u_0=1, u_1, u_2, \ldots, u_n = k), \pi(1), \pi(2), \ldots, \pi(n)},
\end{equation*}
where $\pi(j)$ is supported on $C_j = [u_{j-1} + 1, u_j - 1]$. Let $X$ be a free Poisson stochastic measure. Then for $\mb{X}$ with all $X^{(i)} = X$,
\begin{equation*}
\St_\pi (\mb{X}) = \prod_{j=1}^{n} R_{\pi(j)}(C_j; \mb{X}) \cdot \Delta (B; \mb{X}).
\end{equation*}
\end{Lemma}

\begin{proof}
Let $Z_{i,j}(N) = \sum_{\bar{v} \in ([N]\backslash \set{i})^{\abs{C_j}}_{\pi(j)}} \mb{X}_{(C_j; \bar{v})}$. 
\begin{equation*}
\sum_{i=1}^N X^{(u_0)}_i Z_{i,1} X^{(u_1)}_i Z_{i,2} \ldots Z_{i,n} X^{(u_n)}_i
= \sum_{i=1}^N \sum_{G \subset [n]} (X^{(u_0)}_i, \ldots, X^{(u_n)}_i) \circ Z(G),
\end{equation*}
where $Z(G)$ is a vector of length $n$ such that 
\begin{equation*}
Z(G)_j = 
\begin{cases}
Z_{i,j} - \state{Z_{i,j}}, j \in G, \\
\state{Z_{i,j}}, j \not \in G.
\end{cases}
\end{equation*} 
For $G \neq \emptyset$, at least one of the $Z(G)_j$ is centered. So Lemma \ref{Lem:PoissonCentered} applies and the limit, as $\delta(\mc{S})$ goes to $0$, of the appropriate term is $0$. On the other hand, it follows from \cite[Corollary 2]{Ans00} that for any $i$ the limit of $\state{Z_{i,j}}$ is $R_{\pi_j}(C_j; \mb{X})$. We conclude that, denoting $\bar{\imath} = (i, i, \ldots, i)$, 
\begin{align*}
\lim_{\delta(\mc{S}) \rightarrow 0} \sum_{i=1}^N X^{(u_0)}_i Z_{i,1} X^{(u_1)}_i Z_{i,2} \ldots Z_{i,n} X^{(u_n)}_i 
&= \prod_{j=1}^n R_{\pi(j)}(\mb{X}) \cdot \lim_{\delta(\mc{S}) \rightarrow 0} \sum_{i=1}^N \mb{X}_{(B; \bar{\imath})} \\
&= \prod_{j=1}^n R_{\pi(j)}(\mb{X}) \cdot \Delta(B; \mb{X}),
\end{align*}
where $\Delta(B; \mb{X})$ is well-defined for the free Poisson stochastic measure by \cite[Corollary 4]{Ans00}. On the other hand,
\begin{equation*}
\sum_{i=1}^N X^{(u_0)}_i Z_{i,1} X^{(u_1)}_i Z_{i,2} \ldots Z_{i,n} X^{(u_n)}_i = \sum_{\sigma} \St_\sigma(\mb{X}, \mc{S}),
\end{equation*}
where the sum is taken over all partitions $\sigma$ of $[k]$ which contain the class $B$ and such that for all j, the restriction of $\sigma$ to $C_j$ is $\pi(j)$. The only noncrossing partition satisfying these requirements is $\pi$, so
\begin{equation*}
\lim_{\delta(\mc{S}) \rightarrow 0} \sum_{\sigma} \St_\sigma(\mb{X}, \mc{S}) = \St_\pi(\mb{X}). 
\end{equation*}
\end{proof}

\begin{Notation}
Let $\pi \in \NC(k)$. Then $\pi$ can be written as
\begin{equation*}
\pi = (B_1(\pi), B_2(\pi), \ldots, B_{o(\pi)}(\pi), \mc{I}_1(\pi), \mc{I}_2(\pi), \ldots, \mc{I}_{o(\pi)}(\pi)).
\end{equation*}
Here, $\set{B_i(\pi)}$ are outer classes of $\pi$, listed in increasing order. Denote
\[
\overline{B_i} = \set{j| \exists a,b \in B_i: a \leq j \leq b}
\]
the subset covered by $B_i$, and let $\mc{I}_i(\pi)$ be the restriction of $\pi$ to the set $\overline{B_i} \backslash B_i$ strictly covered by $B_i$. Denote by $\mc{I}_i'(\pi)$ the noncrossing partition $(B_i(\pi), \mc{I}_i(\pi))$. Finally, denote $C(\pi) = [k] \backslash \bigcup_{i=1}^{o(\pi)} B_i$, and let $\mc{I}(\pi)$ be the partition consisting of all inner classes of $\pi$, i.e.\ the restriction of $\pi$ to $C(\pi)$.
\end{Notation}

\begin{Lemma}
\label{lem:prod}
With the above notation, for a consistent $k$-tuple $\mb{X}$ of free stochastic measures, $\Prod_\pi(\mb{X}) = \prod_{i=1}^{o(\pi)} \Prod_{\mc{I}_i'(\pi)}(\mb{X})$.
\end{Lemma}

\begin{proof}
Such a product decomposition is valid for any subdivision $\mc{S}$. 
\end{proof}

\begin{Lemma}
\label{Lem:StPr}
Let $\mb{X}$ be a consistent $k$-tuple of free stochastic measures. Then
\begin{enumerate}
\item
\label{Num:StPr}
The measures $\Prod_\pi(\mb{X})$ and $\St_\pi(\mb{X})$ are related as follows: for $\pi \in \NC(k)$,
\begin{align*}
\Prod_\pi(\mb{X}) &= \sum_{\substack{\sigma \in \NC(k) \\ \sigma \geq \pi}} \St_\sigma(\mb{X}), \\
\St_\pi(\mb{X}) &= \sum_{\substack{\sigma \in \NC(k) \\ \sigma \geq \pi}} \Mob(\pi, \sigma)
\Prod_\sigma(\mb{X}).
\end{align*}
\item
\label{Num:StSt}
Let $\pi_1, \pi_2, \ldots, \pi_n$ be noncrossing partitions such that $\pi = \pi_1 + \pi_2 + \ldots + \pi_n \in \NC(k)$. For each $i$, identify $\pi_i$ with a sub-partition of $\pi$, and let $C_i$ be the support of $\pi_i$ in $[k]$. Denote $\tau \in \NC(n)$ the partition $(C_1, C_2, \ldots, C_n)$. Then
\[
\prod_{i=1}^n \St_{\pi_i}(C_i; \mb{X}) = \sum_{\substack{\sigma \in \NC(k) \\ \sigma \wedge \tau = \pi}} \St_\sigma(\mb{X}).
\]
\end{enumerate}
\end{Lemma}

\begin{proof}
The first statement is based on a purely combinatorial observation that 
\[
\Prod_\pi(\mb{X}, \mc{S}) = \sum_{\substack{\sigma \in \Part(k) \\ \sigma \geq \pi}} \St_\sigma(\mb{X}, \mc{S})
\]
and the fact that $\St_\sigma(\mb{X}) = 0$ for $\sigma \not \in \NC(k)$; see Corollary 1 of \cite{Ans00}. The second statement is based on a purely combinatorial observation that
\[
\prod_{i=1}^n \St_{\pi_i}((C_i; \mb{X}), \mc{S}) = \sum_{\substack{\sigma \in \Part(k) \\ \sigma \wedge \tau = \pi}} \St_\sigma(\mb{X}, \mc{S})
\]
and the same fact.
\end{proof}

\begin{Lemma}
\label{lem:limit}
The limit defining $\Prod_\pi(\mb{X})$ exists in norm if the corresponding limit exists for the free Poisson stochastic measure.
\end{Lemma}

\begin{proof}
Let $\mc{S} = (I_1, I_2, \ldots, I_N)$ be a subdivision of $[0,1)$. Let $\mc{T}$ be another such subdivision, and let $\mc{S} \wedge \mc{T} = (J_1, J_2, \ldots, J_M)$ be their common refinement. Temporarily denote by $p(s)$ the index $i$ such that $J_s \subset I_i$. Denote $A(\mc{S}) = \Prod_\pi(\mb{X}, \mc{S})$, and similarly for $\mc{S} \wedge \mc{T}$.
\begin{align*}
A(\mc{S}) - A(\mc{S} \wedge \mc{T}) 
=& \sum_{\bar{v} \in [N]_{\geq \pi}^k} \sum_{p^{-1}(s_1) = v_1} X^{(1)}_{s_1}(\mc{S} \wedge \mc{T}) \cdots \sum_{p^{-1}(s_k) = v_k} X^{(k)}_{s_k}(\mc{S} \wedge \mc{T}) \\
&\quad- \sum_{\bar{u} \in [M]_{\geq \pi}^k} \mb{X}_{\bar{u}}(\mc{S} \wedge \mc{T}) \\
=& \sum_{\substack{p(\bar{s}) \in [N]_{\geq \pi}^k, \\ \bar{s} \not \in [M]_{\geq \pi}^k}} \mb{X}_{\bar{s}}(\mc{S} \wedge \mc{T}).
\end{align*}
The above expression $A(\mc{S}) - A(\mc{S} \wedge \mc{T})$ is a sum with positive coefficients. Hence so is $((A(\mc{S}) - A(\mc{S} \wedge \mc{T})) (A(\mc{S}) - A(\mc{S} \wedge \mc{T}))^\ast)^n$. Therefore its expectation is a sum over a collection of indices, with weights given by products of $\abs{J_s}$, all of which are independent of the distribution of $\mb{X}$, of free cumulants of $\mb{X}$ of order $2kn$. Each of those free cumulants is bounded in norm by $(16 \norm{\mb{X}})^{2nk}$, where $\norm{\mb{X}} = \max_i \norm{X^{(i)}}$.  Since for the free Poisson process all such cumulants are equal to $1$, the result is at most $(16 \norm{\mb{X}})^{2nk}$ times the corresponding quantity for the free Poisson process, for which we denote $\Prod_\pi(X, \mc{S})$ by $a(\mc{S})$. That is, 
\begin{multline*}
\state{((A(\mc{S}) - A(\mc{S} \wedge \mc{T})) (A(\mc{S}) - A(\mc{S} \wedge \mc{T}))^\ast)^n} \\
\leq (16 \norm{\mb{X}})^{2nk} \state{((a(\mc{S}) - a(\mc{S} \wedge \mc{T})) (a(\mc{S}) - a(\mc{S} \wedge \mc{T}))^\ast)^n},
\end{multline*}
and so
\begin{equation*}
\norm{A(\mc{S}) - A(\mc{S} \wedge \mc{T})}_{2n} \leq (16 \norm{\mb{X}})^k \norm{a(\mc{S}) - a(\mc{S} \wedge \mc{T})}_{2n},
\end{equation*}
which implies in particular that
\begin{equation*}
\norm{A(\mc{S}) - A(\mc{S} \wedge \mc{T})} \leq (16 \norm{\mb{X}})^k \norm{a(\mc{S}) - a(\mc{S} \wedge \mc{T})}.
\end{equation*}
By assumption, the net $a(\mc{S})$ converges in norm, and 
\begin{align*}
\norm{A(\mc{S}) - A(\mc{T})} 
&\leq \norm{A(\mc{S}) - A(\mc{S} \wedge \mc{T})} + \norm{A(\mc{T}) - A(\mc{S} \wedge \mc{T})} \\
&\leq (16 \norm{\mb{X}})^k (\norm{a(\mc{S}) - a(\mc{S} \wedge \mc{T})} + \norm{a(\mc{T}) - a(\mc{S} \wedge \mc{T})}).
\end{align*}
Therefore the net $A(\mc{S})$ is a Cauchy net, and so converges.
\end{proof}

\begin{Cor}
$\Prod_\pi(\mb{X})$, and hence $\St_\pi(\mb{X})$, is well-defined for all $\pi, \mb{X}$.
\end{Cor}

\begin{proof}
Let $X$ be a free Poisson stochastic measure. By Lemma \ref{lem:st1}, $\St_\pi(X)$ is well-defined for $\pi \in \NC(k)$ with a single outer class. Since for such $\pi$ and $\sigma \in \NC(k)$, $\sigma \geq \pi$, $\sigma$ also contains only one outer class, by Lemma \ref{Lem:StPr} part \eqref{Num:StPr} we conclude that for such $\pi$, $\Prod_\pi(X)$ is well-defined as well. By Lemma \ref{lem:prod}, $\Prod_\pi(X)$ is then well-defined for an arbitrary  $\pi \in \NC(k)$, and applying Lemma \ref{Lem:StPr} part \eqref{Num:StPr} again implies that $\St_\pi(X)$ is well-defined for an arbitrary $\pi$ as well. Finally, by Lemma \ref{lem:limit} the same is true for an arbitrary consistent $k$-tuple $\mb{X}$ of free stochastic measures.
\end{proof}

\begin{Cor}
Let $\mb{X}$ be a consistent $k$-tuple of free stochastic measures. For an interval $I$, define $\Delta(\mb{X})(I) = \lim_{\delta(\mc{S}) \rightarrow 0} \St_{\hat{1}}(\mb{X}, \mc{S})$, where $\mc{S}$ is a subdivision of $I$ in place of $[0,1)$. With this notation, $\Delta(\mb{X})$ is a free stochastic measure. 
\end{Cor}

\begin{Lemma}
\label{Lemma:Diagonal}
Let $\mb{X}$ be a consistent $k$-tuple of free stochastic measures. Let $G_1, G_2, \ldots, G_m \subset [k]$, and denote $\mb{X}_G = (X^{(u_1)}, X^{(u_2)}, \ldots X^{(u_{\abs{G}})})$ for $G = (u_1 < u_2 < \ldots < u_{\abs{G}})$. Then the $m$-tuple
\[
(\Delta(\mb{X}_{G_1}), \Delta(\mb{X}_{G_2}), \ldots, \Delta(\mb{X}_{G_m}))
\]
is also consistent.
\end{Lemma}

\begin{proof}
The free increments property and stationarity follow immediately from the corresponding properties of $\mb{X}$. For a general $n$-tuple $\mb{Y}$ of free stochastic measures that has these two properties, by stationarity the continuity property is equivalent to the continuity of the function
\[
t \mapsto \state{Y^{(v_1)}(t) \ldots Y^{(v_l)}(t)}
\]
for all $t, \bar{v}$. By M\"{o}bius inversion, this is equivalent to the continuity of
\[
t \mapsto R(Y^{(v_1)}(t), \ldots,  Y^{(v_l)}(t))
\]
for all $t, \bar{v}$. By additivity and the free increments property, this is equivalent to the continuity of this function, for all $\bar{v}$, at $t=0$, and so to the same property for $M$.

Thus finally, for the $m$-tuple in the hypothesis, it suffices to prove that
\[
\state{\Delta(\mb{X}_{G_1}, t) \Delta(\mb{X}_{G_2}, t) \ldots \Delta(\mb{X}_{G_m}, t)} \rightarrow 0
\]
as $t \rightarrow 0$. Note that we do not need to put in a multi-index $\bar{v}$ since $\set{G_i}_{i=1}^m$ is already an arbitrary collection of subsets of $[k]$. Denote by $\sigma \in \NC(l)$ the partition $(B_1, B_2, \ldots, B_m)$ with interval classes
\begin{equation*}
B_j = \set{(\sum_{s=1}^{j-1} \abs{G_s}) + 1, \ldots, \sum_{s=1}^j \abs{G_s}},
\end{equation*}
and let $\mb{Y} = (\mb{X}_{G_1}, \mb{X}_{G_2}, \ldots, \mb{X}_{G_m})$. Clearly $\mb{Y}$ is a consistent $l$-tuple. Then
\begin{align*}
\state{\Delta(\mb{X}_{G_1}, t) \Delta(\mb{X}_{G_2}, t) \ldots \Delta(\mb{X}_{G_m}, t)}
&= \state{\Prod_\sigma(\mb{Y}(t))} = \sum_{\tau \geq \sigma} \state{\St_\tau(\mb{Y}(t))} \\
&= \sum_{\tau \geq \sigma} R_\tau(\mb{Y}(t)) = \sum_{\tau \geq \sigma} t^{\abs{\tau}} R_\tau(\mb{Y})
\end{align*}
by Lemma \ref{lem:r}, and so goes to $0$ as $t \rightarrow 0$.
\end{proof}

\begin{Lemma}
\label{Lem:Diagonal}
Let $\sigma = (B_1, B_2, \ldots, B_n)$ be an interval partition of $[k]$. Then
\[
\Delta(\Delta(B_1; \mb{X}), \ldots, \Delta(B_n; \mb{X})) = \Delta(\mb{X}).
\]
\end{Lemma}

\begin{proof}
Let $\mc{S} = (I_1, \ldots, I_N)$ be a subdivision of $[0, 1)$. For each $i$, let $\mc{S}_i = (I_{i, 1}, I_{i, 2}, \ldots, I_{i, M_i})$ be a subdivision of $I_i$, and $\mc{T}$ be the subdivision of $[0,1)$ obtained by combining $\set{\mc{S}_i}_{i=1}^N$. Then as $\delta(\mc{S}_1), \ldots, \delta(\mc{S}_N) \rightarrow 0$, also $\delta(\mc{T}) \rightarrow 0$. Therefore 
\[
\lim_{\delta(\mc{S}_1), \ldots, \delta(\mc{S}_N) \rightarrow 0} \Delta(\mb{X}, \mc{T})
= \Delta(\mb{X}),
\]
and so $\Delta(\mb{X})$ is also the limit of the left-hand-side if in addition $\delta(\mc{S}) \rightarrow 0$. Here
\[
\Delta(\mb{X}, \mc{T})
= \sum_{i=1}^N \sum_{s=1}^{M_i} \prod_{t=1}^k X^{(t)}(I_{i, s}).
\]
On the other hand,
\[
\sum_{i=1}^N \prod_{j=1}^n \sum_{s=1}^{M_i} \prod_{t \in B_j} X^{(t)}(I_{i, s})
= \sum_{i=1}^N \prod_{j=1}^n \Delta((B_j; \mb{X}), \mc{S}_i)
\]
and 
\begin{align*}
\lim_{\delta(\mc{S}) \rightarrow 0} \lim_{\delta(\mc{S}_1), \ldots, \delta(\mc{S}_N) \rightarrow 0} \sum_{i=1}^N \prod_{j=1}^n \Delta((B_j; \mb{X}), \mc{S}_i)
&= \lim_{\delta(\mc{S}) \rightarrow 0} \sum_{i=1}^N \prod_{j=1}^n \Delta(B_j; \mb{X}) \\
&= \Delta(\Delta(B_1; \mb{X}), \ldots, \Delta(B_n; \mb{X})).
\end{align*}
Therefore the difference
\[
\Delta(\Delta(B_1; \mb{X}), \ldots, \Delta(B_n; \mb{X})) - \Delta(\mb{X})
\]
is the limit, as $\delta(\mc{S}_1), \ldots, \delta(\mc{S}_N) \rightarrow 0$ and then as $\delta(\mc{S}) \rightarrow 0$, of 
\[
\sum_{i=1}^N \prod_{j=1}^n \sum_{s=1}^{M_i} \prod_{t \in B_j} X^{(t)}(I_{i, s}) - \sum_{i=1}^N \sum_{s=1}^{M_i} \prod_{t=1}^k X^{(t)}(I_{i, s}).
\]
This expression is a sum with positive coefficients. Also, for the free Poisson process,
\[
\Delta(\Delta(B_1; \mb{X}), \ldots, \Delta(B_n; \mb{X})) - \Delta(\mb{X})
= \Delta(X, \ldots, X) - X = 0.
\]
By the same estimates as in Lemma \ref{lem:limit}, the result follows.
\end{proof}

\begin{Lemma}
\label{Lem:Inner}
For $\pi \in \NC(k)$, 
\[
\St_\pi(\mb{X}) = R_{\mc{I}(\pi)}(C(\pi); \mb{X}) \cdot \St_{(B_1(\pi), B_2(\pi), \ldots, B_{o(\pi)}(\pi))}(\bigcup_{i=1}^{o(\pi)} B_i(\pi); \mb{X}).
\]
\end{Lemma}

\begin{proof}
Let $C$ be an inner class of $\pi$, and let $\pi' \in \NC(k - \abs{C})$ be the restriction of $\pi$ to $[k] \backslash C$. Then it suffices to prove that
\[
\St_\pi(\mb{X}) = R(C; \mb{X}) \cdot \St_{\pi'}(([k] \backslash C); \mb{X}).
\] 
Denote $A = \St_\pi(\mb{X}) - R(C; \mb{X}) \cdot \St_{\pi'}(([k] \backslash C); \mb{X})$.
\[
\state{(A A^\ast)^n}
= \sum_{G \subset [2n]} (-R(C; \mb{X}))^{\abs{G}} \state{\St_{\pi_1}(\mb{X}_1) \St_{\pi_2}(\mb{X}_2) \ldots \St_{\pi_{2n}}(\mb{X}_{2n})},
\]
where
\begin{align*}
\text{If } & j \not \in G, j \text{ odd, then } \pi_j = \pi, \mb{X}_j = \mb{X}. \\
\text{If } & j \not \in G, j \text{ even, then } \pi_j = \pi^{op}, \mb{X}_j =  \mb{X}^{op}. \\
\text{If } & j \in G, j \text{ odd, then } \pi_j = \pi', \mb{X}_j = (([k] \backslash C); \mb{X}). \\
\text{If } & j \in G, j \text{ even, then } \pi_j = (\pi')^{op}, \mb{X}_j = (([k] \backslash C); \mb{X})^{op}.
\end{align*}
Denote $\pi_G = \pi_1 + \pi_2 + \ldots + \pi_{2n}$. Let $C_i(G)$ be the support of $\pi_i$ identified as a sub-partition of $\pi_G$, and let $\tau_G = (C_1(G), C_2(G), \ldots, C_{2n}(G))$. Then by part \eqref{Num:StSt} of Lemma \ref{Lem:StPr},
\begin{align*}
\norm{A}_{2n}^{2n}
&= \sum_{G \subset [2n]} (-1)^{\abs{G}} R(C; \mb{X})^{\abs{G}} \sum_{\substack{\sigma \in \NC(2n k - \abs{G} \cdot \abs{C}) \\ \sigma \wedge \tau_G = \pi_G}}  \state{\St_\sigma(\mb{X}_1, \mb{X}_2, \ldots, \mb{X}_{2n})} \\
&= \sum_{G \subset [2n]} (-1)^{\abs{G}} \sum_{\substack{\sigma \in \NC(2n k - \abs{G} \cdot \abs{C}) \\ \sigma \wedge \tau_G = \pi_G}} R(C; \mb{X})^{\abs{G}} R_\sigma(\mb{X}_1, \mb{X}_2, \ldots, \mb{X}_{2n}).
\end{align*}
Fix $G \subset [2n]$. Let $\sigma \in \NC(2nk)$, $\sigma \wedge \tau_\emptyset = \pi_\emptyset$, where $\tau_\emptyset = \hat{1}_k + \ldots + \hat{1}_k$ and $\pi_\emptyset = \pi + \pi^{op} + \pi + \ldots + \pi^{op}$. Denote $C^{op} = (k+1-C)$ the class of $\pi^{op}$ corresponding to $C$. Since $C$ is an inner class of $\pi$, the condition $\sigma \wedge \tau_\emptyset = \pi_\emptyset$ implies that $(2j k + C)$ and $((2j+1)k + C^{op})$ are classes of $\sigma$ for $0 \leq j < n$. Let $g_G$ map such a $\sigma$ to the partition in $\NC(2n k - \abs{G} \cdot \abs{C})$ obtained by removing from $\sigma$ the classes $(2j k + C)$ for $(2j+1) \in G$ and $((2j+1)k + C^{op})$ for $(2j+2) \in G$. It is easy to see that $g_G$ is a bijection onto $\set{\sigma \in \NC(2n k - \abs{G} \cdot \abs{C}) | \sigma \wedge \tau_G = \pi_G}$, and that $R(C; \mb{X})^{\abs{G}} R_{g_G(\sigma)}(\mb{X}_1, \mb{X}_2, \ldots, \mb{X}_{2n}) = R_\sigma(\mb{X}, \mb{X}, \ldots, \mb{X})$. Therefore
\[
\norm{A}_{2n}^{2n}
= \sum_{G \subset [2n]} (-1)^{\abs{G}} \sum_{\substack{\sigma \in \NC(2n k) \\ \sigma \wedge \tau_\emptyset = \pi_\emptyset}} R_\sigma(\mb{X}, \mb{X}, \ldots, \mb{X}) = 0
\]
since the first sum equals to $0$.
\end{proof}

\begin{proof}[Proof of the Main Theorem]
The statement of the theorem holds for $\pi = \hat{0}_k$. From now on, assume $\pi > \hat{0}_k$. The proof will proceed by induction on $k$. The statement of the theorem is vacuous for $k=1$; assume that it holds for all tuples of less than $k$ elements.

By Lemma \ref{Lem:Inner},
\[
\St_{\mc{I}_i'(\pi)} (\overline{B_i(\pi)}; \mb{X}) 
= R_{\mc{I}_i(\pi)}(\overline{B_i(\pi)} \backslash B_i(\pi); \mb{X}) \cdot \Delta(B_i(\pi); \mb{X}).
\]
Therefore 
\begin{align*}
\Prod_{\mc{I}_i'(\pi)} (\overline{B_i(\pi)}; \mb{X}) 
&= \sum_{\sigma_i \geq \mc{I}_i'(\pi)} \St_{\sigma_i} (\overline{B_i(\pi)}; \mb{X}) \\
&= \sum_{\sigma_i \geq \mc{I}_i'(\pi)} (R_{\mc{I}(\sigma_i)}(\overline{B(\sigma_i)} \backslash B(\sigma_i); \mb{X}) \cdot \Delta(B(\sigma_i); \mb{X})).
\end{align*}
Then by Lemma \ref{lem:prod},
\begin{align*}
\Prod_\pi(\mb{X}) 
&= \prod_{i=1}^{o(\pi)} \Prod_{\mc{I}_i'(\pi)} (\overline{B_i(\pi)}; \mb{X}) \\
&= \prod_{i=1}^{o(\pi)} \sum_{\sigma_i \geq \mc{I}_i'(\pi)} (R_{\mc{I}(\sigma_i)}(\overline{B(\sigma_i)} \backslash B(\sigma_i); \mb{X}) \cdot \Delta(B(\sigma_i); \mb{X})) \\
&= \sum_{\substack{\sigma \geq \pi \\ \forall i: \overline{B_i(\sigma)} = \overline{B_i(\pi)}}} R_{\mc{I}(\sigma)}(C(\sigma); \mb{X}) \prod_{j=1}^{o(\pi)} \Delta(B_j(\sigma); \mb{X}) \\
&= \sum_{\substack{\sigma \geq \pi \\ \forall i: \overline{B_i(\sigma)} = \overline{B_i(\pi)}}} R_{\mc{I}(\sigma)}(C(\sigma); \mb{X}) \\
&\qquad \times \Prod_{\hat{0}_{o(\pi)}}(\Delta(B_1(\sigma); \mb{X}), \ldots, \Delta(B_{o(\pi)}(\sigma); \mb{X})).
\end{align*}
In its turn,
\begin{multline*}
\Prod_{\hat{0}_{o(\pi)}}(\Delta(B_1(\sigma); \mb{X}), \ldots, \Delta(B_{o(\pi)}(\sigma); \mb{X})) \\
= \sum_{\rho \in \NC(o(\pi))} \St_\rho(\Delta(B_1(\sigma); \mb{X}), \ldots, \Delta(B_{o(\pi)}(\sigma); \mb{X})).
\end{multline*}
Since $\pi > \hat{0}_k$, $\sigma$ has at most $k-1$ classes, so the induction hypothesis applies to $\mb{Y} = (\Delta(B_1(\sigma); \mb{X}), \ldots, \Delta(B_{o(\sigma)}(\sigma); \mb{X}))$. Thus 
\begin{equation}
\label{strho}
\St_\rho(\mb{Y})
= R_{\mc{I}(\rho)}(C(\rho); \mb{Y}) \cdot \psi(\Delta(B_1(\rho); \mb{Y}), \ldots, \Delta(B_{o(\rho)}(\rho); \mb{Y})).
\end{equation}
Define the map $f: \NC(o(\pi)) \times \set{\sigma \in \NC(k)| \sigma \geq \pi, \forall i: \overline{B_i(\sigma)} = \overline{B_i(\pi)}} \rightarrow \NC(k)$ by $i \stackrel{f(\rho, \sigma)}{\sim} j \Leftrightarrow ((i \stackrel{\sigma}{\sim} j) \text{ or } (i \in B_s(\sigma), j \in B_t(\sigma), s \stackrel{\rho}{\sim} t))$. Note that the outer classes of $f(\rho, \sigma)$ are in one-to-one correspondence with the outer classes of $\rho$, and each inner class of $f(\rho, \sigma)$ corresponds to a unique inner class of either $\rho$ or $\sigma$. It is easy to see that $f$ is in fact a bijection onto $\set{\tau \in \NC(k)| \tau \geq \pi}$. Combining equation \eqref{strho} with Lemma \ref{Lem:Diagonal}, we see that
\begin{multline*}
R_{\mc{I}(\sigma)}(C(\sigma); \mb{X}) \cdot \St_\rho(\Delta(B_1(\sigma); \mb{X}), \ldots, \Delta(B_{o(\pi)}(\sigma); \mb{X})) \\
= R_{\mc{I}(\tau)}(C(\tau); \mb{X}) \cdot \psi(\Delta(B_1(\tau), \mb{X}), \ldots, \Delta(B_{o(\tau)}(\tau); \mb{X})),
\end{multline*}
with $\tau = f(\rho, \sigma)$. Therefore
\begin{equation*}
\Prod_\pi(\mb{X}) 
= \sum_{\tau \geq \pi} R_{\mc{I}(\tau)}(C(\tau); \mb{X})
\times \psi(\Delta(B_1(\tau), \mb{X}), \Delta(B_2(\tau); \mb{X}), \ldots, \Delta(B_{o(\tau)}(\tau); \mb{X})).
\end{equation*}
On the other hand, for all $\pi$, $\Prod_\pi(\mb{X}) = \sum_{\tau \geq \pi} \St_\pi(\mb{X})$. Note that the M\"{o}bius inversion formula for $\pi > \hat{0}_k$ involves only $\sigma > \hat{0}_k$. Therefore, applying this formula,
\begin{align*}
\St_\pi (\mb{X})
&= R_{\mc{I}(\pi)}(C(\pi); \mb{X}) \times \psi(\Delta(B_1(\pi), \mb{X}), \Delta(B_2(\pi); \mb{X}), \ldots, \Delta(B_{o(\pi)}(\pi); \mb{X})) \\
&= \prod_{i=1}^{i(\pi)} R(C_i; \mb{X}) \cdot \psi(\Delta(B_1; \mb{X}), \Delta(B_2; \mb{X}), \ldots, \Delta(B_{o(\pi)}; \mb{X})). 
\end{align*}
\end{proof}

\begin{Remark}[Higher-dimensional analogs]
\label{Remark:hda}
The Main Theorem gives a complete description of the higher stochastic measures $\St_\pi$ as given in Definition \ref{Defn:fsm}. However, under the original definitions of \cite{RW97} (modified for processes with freely independent increments) these only correspond to values on cubes, hence their dependence on only $1$ and not $k$ parameters. In this remark we briefly describe how one could extend the definition to more general rectangles of the form $\mb{I} = [a_1, b_1) \times [a_2, b_2) \times \cdots \times [a_k, b_k)$. It is clear that it suffices to give the definition only for the case when for $1 \leq i, j \leq k$ the intervals $[a_i, b_i)$ and $[a_j, b_j)$ are either disjoint or the same (one then needs to show that the resulting definition is consistent). Assume that the rectangle $\mb{I}$ is of this form. Then we can define a partition $\pi(\mb{I}) \in \Part(k)$ by $i \stackrel{\pi(\mb{I})}{\sim} j \Leftrightarrow [a_i, b_i) = [a_j, b_j)$. Let $\pi(\mb{I})$ have classes $B_1, B_2, \ldots, B_l$. Let $c(i)$ be the index such that $i \in B_{c(i)}, 1 \leq i \leq k$. Let $X$ be a free stochastic measure, and $\mb{X}$ a $k$-tuple of free stochastic measures given by $X^{(j)}([a,b)) = X([a - a_j, b - a_j))$. The conditions on $a_i, b_i$ imply that this $k$-tuple is consistent. Let $\mb{S} = \set{\mc{S}_j}$ be subdivisions of $[a_{c^{-1}(j)}, b_{c^{-1}(j)}), 1 \leq j \leq l$ into intervals $I_{j, s}$. For $\sigma \in \NC(k)$, denote $\mb{S}_\sigma = \set{\bar{v} \in \mf{N}^k: (\prod_{i=1}^k I_{c(i), v_i}) \cap \mf{R}^k_\sigma \neq \emptyset}$. Note that if $\pi(\mb{I}) = \hat{1}$ and $\mc{S}$ is a single subdivision with $N$ classes, $\mc{S}_\sigma = [N]^k_\sigma$. Define  
\begin{equation*}
\St_\sigma(X, \mb{S}) = \sum_{\bar{v} \in \mb{S}_\sigma} \prod_{i=1}^k X(I_{i, v_i})
\end{equation*} 
and $\St(\mb{I}) = \lim_{\delta(\mb{S}) \rightarrow 0} \St_\sigma(X, \mb{S})$. It follows immediately that $\St(\mb{I}) = 0$ unless $\sigma \leq \pi(\mb{I})$. Indeed, if $\sigma \not \leq \pi(\mb{I})$ then for any subdivision $\mb{S}$, $\mb{S}_\sigma = \emptyset$. 
\end{Remark}

\end{document}